\documentclass[a4paper]{amsart}
\usepackage[latin1]{inputenc}
\usepackage{amsmath}
\usepackage{amsthm}
\usepackage{latexsym}
\usepackage{amssymb}
\usepackage{wasysym}
\usepackage{setspace}
\usepackage{mathrsfs}

\newtheorem{definition}{Definition}[section]
\newtheorem{theorem}[definition]{Theorem}
\newtheorem{corollary}[definition]{Corollary}

\newtheorem{lemma}[definition]{Lemma}
   
\allowdisplaybreaks
\numberwithin{equation}{section}

\begin{document}
\title{Definitions of Almost Periodicity}
\author{H. Günzler}
\begin{abstract}
Various versions of the classical definitions of (one- and twosided) almost periodicity for functions on groups with values in a uniform space are formulated and their equivalence is shown.
\end{abstract}
\maketitle
\section{Introduction}
The first generalization of continuous periodic $f: \mathbb{R} \rightarrow \mathbb{C}$ in the direction of almost periodicity seems due to Bohl 1893 and independently Esclangon 1904, see \cite{lit13}, \cite{lit29}; these "quasi-periodic" functions, which appeared in astronomical perturbation problems, are in their simplest form $\sum_{j=1}^n c_j \sin (a_j t) $ with given $a=(a_1, \dots, a_n) \in \mathbb{R}^n$, in genereal they are, with given fixed $a \in \mathbb{R}^n$, those that can be approximated, uniformly on $\mathbb{R}$, by expressions $\sum_{k=1}^n c_k A_k$ with $c_k \in \mathbb{C}$ and $A_k(t)=e^{i(m_1 a_1 + \dots{} + m_n a_n)t}$ with $m_j$ integers (or restrictions of continuous functions on a $n$-dimensionable torus) \cite[p. 145/146]{lit38}, \cite{lit29}.

Bohr's almost periodic (=ap) functions of 1923 \cite{lit14}, \cite{lit15}, \cite{lit16} were defined as those continuous $f: \mathbb{R} \rightarrow \mathbb{C}$ which had, for any $\epsilon > 0 $, sufficiently many $\epsilon$-periods $\tau$, i.e. $|f(t+\tau)-f(t)|< \epsilon$ for all $t\in \mathbb{R}$ (definitions 3.11 and 3.7 below).

Bochner 1927 characterized these ap functions as those $ f: \mathbb{R} \rightarrow \mathbb{C}$ which for any sequence $(s_m)$ from $\mathbb{R}$ admitted a subsequence $(s_{m_n})$ so that $ (f(s_{m_n} + \cdot \ ))$ converges uniformly on $\mathbb{R}$ \cite{lit9}; in \cite{lit10} he extended Bohr's theory to $f: \mathbb{R} \rightarrow \text{Banach space}$.

With the Haar measure Bohr's theory has been extended to locally compact groups and then 1934 by von Neumann \cite{lit42} to $\mathbb{C}$-valued functions on arbitrary groups. Maak's definition is a combinatorial reformulation of von Neumann's definition, this simplifies the proofs of the main theorems (existence of a mean, Fourier series, approximation theorem) \cite{lit36}.

Finally in \cite{lit11} Bochner and von Neumann extended all this to $f: G \rightarrow X$ with arbitrary group $G$ and $X=$ locally convex topologically complete linear space (exactly for such an extension von Neumann introduced these $X$ \cite{lit43}). Formally subsumed by the Bochner von Neumann theory but introduced independently are the weakly ap functions of Amerio \cite{lit1} ($f:\mathbb{R} \rightarrow B$ Banach space with $f$ ap in the weak topology), needed for his proof that the solutions of the mixed problem for the wave equation are time-ap \cite{lit2}. (These functions should not be confused with the weakly ap (complex valued) functions in the sense of Eberlein, see e.g. \cite{lit20}, \cite{lit8}.)

Though the definitions of ap of Bohr, Bochner, von Neumann/Maak are quite different, implicit in the papers mentioned above is their equivalence.

As documented by several recent articles, this seems no longer general knowledge; also, these equivalences have not been formulated explicitly, especially with respect to what range spaces $X$ are admissible; for non-abelian groups, even von Neumann was not aware that his left ap implies right ap (this is due to Maak \cite[p. 153]{lit38}).

These were some of the reasons for writing this note. Also, in the non-abelian case, some of the characterizations in Theorem 4.1 below seem to be new.\\

\emph{Not} touched upon this note:\\
A) Characterization via uniform approximation by "trigonometric polynomials" (with coefficients of unitary irreducible finite-dimensional complex representations of $G$): For this the existence of an (ergodic) mean is needed, meaning restrictive conditions on the range space $X$ \cite{lit43}, \cite{lit11}, \cite{lit37}, \cite[p. 85]{lit22}, \cite{lit23}, \cite[p. 198 Corollary]{lit26}, \cite[p. 112]{lit27}\\

\noindent B) For further (also "non-uniform") characterizations of ap see \cite[Theorem 2]{lit12}, \cite{lit41}, \cite{lit45}.\\

\noindent C) Some of the results below hold also for semigroups, with a sharpening of Maak's definition of ap, to get still an approximation theorem: \cite[p. 45/46]{lit39}, \cite{lit40}, \cite[Anhang]{lit22}, \cite{lit23}, \cite[p. 198]{lit26}, \cite[p. 112]{lit27}. \\

\noindent D) There is an extensive literature on "almost periodic" functions on semigroups, there von Neumann's definition is used verbatim for semigroups, getting only some kind of asymptotic almost periodicity, they are no longer almost periodic in the intuitive sense as in the group case: \cite{lit20}, \cite{lit8} for example. Already for $f: [0, \infty) \rightarrow \mathbb{C}$ Bohr's definition gives $f=g\big|_{[0,\infty)}$ with unique $g\in AP(\mathbb{R},\mathbb{C})$, but Bochner's definition gives asymptotic ap $f=g\big|_{[0,\infty)} + h$, $h$ continuous with $h(t) \rightarrow 0$ as $t \rightarrow \infty$ \cite[p. 96 footnote $^{11}$ Satz]{lit24}.\\

\noindent E) For $G=\mathbb{R}^n$ (or locally compact $G$) there are several "concrete" gereralizations of Bohr's ap functions which include discontinuous periodic functions: Stepanoff-, Weyl-, Besicovitch-ap functions, ap distributions, for these see e.g. \cite[p. 228-231]{lit38} and the references there, \cite[Chapt VI, § 9 p. 206-208]{lit46}, \cite[(1.9), (2.17), (2.19), Example 3.2]{lit7}, also the references in \cite{lit3}. E.g. by \cite[(3.6), Prop. 3.8, Remark 3.11]{lit5}, \cite[(1.7), (2.19), Ex. 3.2, (3.3)]{lit6} one has, if $1<p<\infty$, with all "$\subset$" strict
\begin{multline}
AP \subset S^pAP  \subset S^1AP  \subset 1.\text{mean class } MAP  \subset M^2AP  \subset \\ \dots  \subset \bigcup_{n=1}^{\infty} M^nAP =
				\mathscr{D}^{'}_{AP} \cap L_{loc}^1 \subset \mathscr{D}^{'}_{AP} = \mathscr{D}^{'}_{S^pAP}=\mathscr{D}^{'}_{S^1AP}  \subset \mathscr{S}'.
\end{multline}

\noindent F) "Abstract" generalizations of ap are almost automorphic, Levitan-ap, recurrent, distal, minimal, Poisson stable, W-ap, $\omega$-ap to mention some, see e.g. Levin 1948 \cite{lit33}, Levitan 1953 \cite{lit34}, Maak 1950 \cite[p. 231 $^{39}$]{lit38}, Maak 1952 \cite[p. 43]{lit39}, Bochner 1962 \cite{lit12}, Reich 1970 \cite{lit44}, Reich 1977 \cite{lit45}, Levitan-Zhikov 1982 \cite[p. 33, 53/54, 63, 64, 76, 80]{lit35}, Milnes 1986 \cite{lit41}, Berglund-Junghenn-Milnes 1989 \cite[Chapt. IV]{lit8}, Basit-Günzler 2005 \cite[p. 427, 429, 430 (3.3)]{lit6}.\\

\noindent G) Also many perturbations of ap functions have been considered, i.e. AP+$N$ with $N$ some nullspace, e.g. Fréchet's asymptotic ap functions \cite[p. 35]{lit48}, \cite[p 43]{lit39}, \cite[p. 30 Theorem 2.22]{lit20}, \cite[p. 157 Theorem 3.13, Example 3.14]{lit8}, \cite[p. 425 and Examples 5.4-5.9]{lit6}, subsuming Zhang's pseudo ap functions, and the references in these books/papers. (See also D)).
\section{Assumptions and Notation}
In the following $\mathbb{N}:=\{1,2,\dots,\}$ and $\mathbb{R} \text{ resp. } \mathbb{C}$ denotes the real resp. complex field, in the usual topology. Topological notions are as in Kelley \cite{lit31}.

$G$ always denotes  a topological group, i.e. a group and a topological space with $(x,y) \mapsto xy$ resp. $x \mapsto x^{-1}$ $(x,y)$-continuous on $G \times G$ resp. $x$-continuous on $G$ (e.g. \cite[p 16]{lit30}); no separation, countability or completeness conditions are assumed; the discrete topology and thus arbitrary abstract groups are included. $\mathscr{N}=\mathscr{N} (e)$ will denote the set of neighborhoods of the unit $e$ of the topology of $G$.

$X$ or $(X,\mathscr{U})$ denotes a nonempty uniform space \cite[p. 176]{lit31} with set of vicinities $\mathscr{U}$, and unique topology on $X$ induced by $\mathscr{U}$ \cite[p. 178]{lit31}; again no separation, countability or completeness axioms are assumed. If $\emptyset \ne M \subset X$, $M$ with the restriction $\mathscr{U}\big|_{M}= \{U\cap (M\times M) \mid U\in \mathscr{U}\}$ is again a uniform space.

Standard assumption in all of the following will be
\begin{equation}
G \text{  is a topological group, } X \text{ a uniform space with uniformity  } \mathscr{U}. \end{equation}
Later we will repeatedly use \cite[p. 176, 179]{lit31}.
\begin{equation} \text{To } U \in \mathscr{U} \text{ exists }  V \in \mathscr{U} \text{ with } V^{-1}=V \subset V \circ V \subset U.\end{equation}
\begin{equation} Y:= C(G,X)= \{f: G \rightarrow X \mid f \text{ continuous in the } \mathscr{U} \text{-topology} \};\end{equation}
in $Y$ uniform convergence is defined with
\begin{equation} \mathscr{V}:= \{ V \subset Y \times Y \mid \text{ to } V \text{ exists } U\in \mathscr{U} \text{ with } V_U \subset V\},\end{equation}
\begin{equation} V_U:= \{(f,g) \in Y \times Y \mid (f(t), g(t)) \in U \text{ for all } t \in G\};
\end{equation}
$(Y,\mathscr{V})$ is again a uniform space, and so is $C(G\times G, X)$ with corresponding $\mathscr{V}_{\times}$.\\
We say that an $M \subset X$ satisfies \emph{von Neumann's countability axiom} $(A_0)$ \cite[p. 4 Def. 2b (2)]{lit43} if the following holds\\
	\begin{align*}(A_0) \ \text{ there exist  } U_n \in \mathscr{U}, n\in \mathbb{N} \text{ with } (M \times M) \cap \left(\bigcap_{n=1}^{\infty} U_n\right) = (\bigcap_{U\in \mathscr{U}}U) \cap (M\times M)\end{align*} 
	
	\begin{flalign*} (A_1) &&  &\text{  means the existence of } U_n \in \mathscr{U}, n\in \mathbb{N}, \text{  so that for each } U \in \mathscr{U} \\&& &\text{ the } U\cap(M\times M) \text{ contains some } U_n \cap (M\times M).\end{flalign*}
An $M \subset X$ is called \emph{totally bounded} if to each $U \in \mathscr{U}$ there exist $n\in\mathbb{N},$$\ x_1, \dots, x_n \in M$ such that:
\begin{equation}M\subset \bigcup_{j=1}^n U(x_j), \ U(x):= \{y\in X \mid (x,y) \in U\} \text{ \cite[p. 198]{lit31}.}
\end{equation}
A uniform $X$ is called \emph{topologically complete} if each closed totally bounded set $M$ from X is complete (each Cauchy net from $M$ converges to some $x \in M$).\\
\emph{Translates} of an $f:G \rightarrow X$ are defined by
\begin{equation}
f_a(t):= f(ta), \ {}_af(t):=f(at), \  a, t \in G.
\end{equation}
ap means almost periodic, $f(M):=\{f(t) \mid t\in M\}$.\\
$A\Rightarrow B$ resp. $A \Leftrightarrow B$ mean $A$ implies $B$ resp. $A$ is equivalent with $B$.\\

Later we will need
\begin{theorem}
For a subset $M$ of a uniform space $(X, \mathscr{U})$ the following 3 conditions are equivalent:
\begin{enumerate}
\item[(a)] $M$ is totally bounded.
\item[(b)] To each net $(x_i)_{i\in I}$ with all $x_i \in M$ there exists a Cauchy subnet $(x_{i_j})_{j\in J}$.
\item[(c)] To each sequence $(x_n)_{n \in \mathbb{N}}$ with all $x_n\in M$ there exists a Cauchy subnet $(x_{n_j})_{j\in J}$.
\end{enumerate}
\end{theorem}
Here a net $(y_j)_{j\in J}$ from $(X, \mathscr{U})$ is called a \emph{Cauchy net} if to each $U\in \mathscr{U}$ there exists a $j_U\in J$ with $(y_k, y_l) \in U$ for all $k, l \in J$ with $j_U \leq k, j_U\leq l$.
\begin{proof}
See Kelley \cite{lit31}, proof of 32 Theorem, p. 198/199.
\end{proof}

\section{Definitions of Almost Periodicity}
In the following, $G$ is a topological group, $X=(X, \mathscr{U})$ a uniform space, $Y=C(G,X)$ of (2.3).
\begin{definition}
$f:G\rightarrow X$ is called 
\begin{enumerate}
\item[a)] \emph{left Maak almost periodic (LMap)} if $f\in C(G, X)$ and to each $U\in \mathscr{U}$ there exist $n\in\mathbb{N}$ and $A_1, \dots, A_n \subset G$ with $\bigcup_{j=1}^n A_j=G$ so that if $u,v \in A_k$ for some $k$ one has $(f(ut), f(vt)) \in U$ for all $t \in G$. (The $A_j$ need not be disjoint.)
\item[b)]right Maak almost periodic \emph{(RMap)} if $f \in C(G,X)$ and to each $U\in \mathscr{U}$ there are $n\in \mathbb{N}, \ A_1, \dots, A_n \subset G$ with $ \bigcup_{j=1}^n A_j=G$ so that if $u,v\in A_k$ for some $k$ one has $(f(tu),f(tv))\in U$ for all $t\in G$.
\item[c)] \emph{LRMap}  if $f$ is LMap and RMap.
\item[d)] Maak almost periodic \emph{(Map)} if $f \in C(G,X)$ and to each $U\in \mathscr{U}$ exist $n\in\mathbb{N},\ A_1, \dots, A_n \subset G$ with $\bigcup_{j=1}^n A_j =G$, so that if $u,v \in A_k$ for some $k$, $(f(sut),f(svt))\in U$ for all $s,t \in G$.
\end{enumerate}
\end{definition}\setcounter{definition}{1}
\begin{definition}\label{def3.2a}
$f:G\rightarrow X$ is called 
\begin{enumerate}
\item[a)] \emph{left Neumann almost periodic (LNap)} if $f\in C(G,X)$ and the set ${}_G f:= \{{}_sf \mid s\in G\}$ of left translates of $f$ is totally bounded in the uniform space $C(G,X)$.
\item[b)]  right Neumann ap \emph{(RNap)} if $f\in C(G,X)$ and the set $f_G:=\{f_s \mid s\in G\}$ of right translates is totally bounded in $C(G,X)$.
\item[c)] \emph{LRNap} if $f$ is LNap and RNap.
\item[d)] Neumann almost periodic \emph{(Nap)} if the set $\{f_{[u]} \mid u\in \mathscr{U}\}$ is totally bounded in the uniform space $C(G\times G, X)$, where  $f_{[u]}(s,t):= f(sut), \ s,t \in G$ and the uniformity for $C(G\times G,X)$ is defined as in section 2 (uniform convergence on $G\times G$).
\end{enumerate}
\end{definition}

The "ap" von Neumann used in \cite[Def. 1, p. 448]{lit42} is our LRNap of Definition 3.2(c).

The Nap used here should not be confused with the Nap introduced by Levitan (see \cite{lit33}, \cite{lit34}, \cite[p. 53/54]{lit35}), these more general functions are now called Levitan almost periodic (e.g. \cite[p. 219/220]{lit44}).

\begin{definition}
$f:G \rightarrow X$ is called
\begin{enumerate}
\item[a)] \emph{left net-net almost periodic (Lnnap)} if $f\in C(G,X)$ and each net $({}_{s_j} f)_{j\in J}$ of left translates of $f, s_j \in G$ has a Cauchy subnet in $C(G,X)$.
\item[b)] \emph{Rnnap} if $f \in Y$ and each net $(f_{s_j})_{j\in J}$ has a Cauchy subnet in $Y=C(G,X)$.
\item[c)] \emph{LRnnap} if $f$ is Lnnap and Rnnap.
\item[d)] \emph{nnap} if $f\in C(G,X)$ and each net $(f_{[s_j]})_{j\in J}$ has a Cauchy subnet in \mbox{$C(G\times G,X)$} (see Definition 3.2 (d)).
\end{enumerate}
\end{definition}

Cauchy net is defined after Theorem 2.1.

\begin{definition}
$f:G\rightarrow X$ is called
\begin{enumerate}
\item[a)] \emph{left sequence-net almost periodic (Lsnap)} means $f\in C(G,X)$ and each sequence $({}_{s_n} f)_{n\in \mathbb{N}}$ has a Cauchy subnet in $C(G,X)$.
\end{enumerate}
\end{definition}
\noindent \textbf{Definitions 3.4 b), c), d)} for \emph{Rsnap, LRsnap, snap} are analogous to Definitions 3.3b, c, d.
 
\begin{definition}
$f: G\rightarrow X$ is called
\begin{enumerate}
\item[a)] \emph{left sequence-sequence almost periodic (Lssap)} if $f\in C(G,X)$ and each sequence $({}_{s_n} f)_{n\in \mathbb{N}}$ has a Cauchy subsequence $({}_{s_{n_k}}f)_{k\in \mathbb{N}}$ in $C(G,X)$.
\end{enumerate}
\end{definition}
This, for abelian $G$, is Bochner's Definition \cite[p. 143, Satz XXI, $f$ "normal"]{lit9}, \cite[p.154]{lit10}.

Since the "B" is needed later for Bohr's definition, we use nn, sn, ss.\\

\noindent \textbf{Definitions 3.5 b), c), d)} are again analogous to Definitions 3.3b, c, d.

\begin{definition}
$M\subset G$ is called
\begin{enumerate}
\item[a)]  \emph{left F-relatively dense (LFrd)} if there exists a finite $F\subset G$ with $FM:=\{uv \mid u\in F, v \in M\}=G$.
\item[b)]  \emph{right F-relatively dense (RFrd)} if there exists a finite $F\subset G$ with $MF=G$.
\item[c)]  \emph{LRFrd} if it is LFrd and RFrd.
\item[d)]  \emph{weakly F-relatively dense} if there exist two finite $F_1, F_2$ with $F_1MF_2=G$.
\end{enumerate}
\end{definition}

\begin{definition}
$M \subset G$ is called
\begin{enumerate}
\item[a)] \emph{left relatively dense (Lrd)} if there exists a compact $K\subset G$ with $KM=G$.
\item[b)] \emph{Rrd} if $MK=G$ for some compact $K\subset G$.
\item[c)] \emph{LRrd} if $M$ is Lrd and Rrd.
\item[d)]  
\begin{sloppypar}
	\emph{weakly relatively dense (weakly rd)} if there are compact $K_1, K_2$ with $K_1MK_2=G$.
\end{sloppypar}

\end{enumerate}
\end{definition}

For a still sufficient weakening see Bohr \cite{lit19}.
\begin{definition}\label{def3.8}
For $f\in C(G,X), U\in \mathscr{U}$,
	\begin{align*}
	P_L(f,U)&:=\{\tau \in G \mid (f(\tau x), f(x))\in U \text{ for all } x \in G\}\\
	P_R(f,U)&:=\{\tau \in G \mid (f( x \tau), f(x))\in U \text{ for all } x \in G\}\\	
	P(f,U)&:=\{\tau \in G \mid (f(x\tau y), f(xy))\in U \text{ for all } x,y \in G\}
	\end{align*}
\end{definition}
\begin{definition}\label{def3.9a}
$f:G\rightarrow X$ is called
\begin{enumerate}
\item[a)]  \emph{left F-Bohr almost periodic (LFBap)} if $f\in C(G,X)$ and for each $U\in \mathscr{U}$ the set $P_L(f,U)$ of left U-periods is left F-relatively dense.
\item[b)] \emph{RFBap} if $f\in C(G,X)$ and $P_R(f,U)$ is RFrd for each $U\in \mathscr{U}$. 
\item[c)] \emph{LRFBap} if f is LFBap and RFBap.
\item[d)]  \emph{FBap} if $P(f,U)$ is weakly Frd for all $U\in \mathscr{U}$.
\end{enumerate}
\end{definition}

\begin{lemma}\label{lemma3.10}
For $G,X, \mathscr{U}$ as in (2.1) and any $f:G\rightarrow X$, the following two statements are equivalent:
	\begin{enumerate}
	\item[(3.1)] For any $U \in \mathscr{U}$, the $P_L(f,U)$ of Definition 3.8 is left F-relatively dense (Definition 	3.6a)
	\item[(3.2)] For any $U \in \mathscr{U}$ the $P_L(f,U)$ is right F-relatively dense.  
	\end{enumerate}
The same holds for the $P_R$ and also for the $P$.
\begin{proof}
"(3.1) $\Rightarrow$ (3.2)": If with (2.2) $V=V^{-1} \in \mathscr{U}, P_L:=P_L(f,V)$, then $\tau \in P_L$ implies $(f(\tau^{-1}s), f(s)) \in V^{-1}=V$ or $ \tau^{-1}\in P_L$, so $G=FP_L$ gives $$G=G^{-1}=P_L^{-1} F^{-1} \subset F^{-1}P_L\subset G,$$ so (3.2) holds.

"(3.1) $\Leftarrow$ (3.2)" and the $P_r$ and $P$ cases follow similarly.
\end{proof}
\end{lemma}

So if in Definition 3.9a the "$P_L$ left Frd" is replaced by "$P_L$ right Frd", one gets the same function class; similarly for Definition 3.9b.
\begin{definition}\label{def3.11a}
$f:G\rightarrow X$ is called
\begin{enumerate}
\item[a)]  \emph{left Bohr almost periodic (LBap)} if $f\in C(G,X)$ and for each $U\in \mathscr{U}$ the $P_L(f,U)$ is left-relatively dense.
\item[b)] \emph{RBap} if $f\in C(G,X)$ and $P_R(f,U)$ is Rrd for each $U\in\mathscr{U}$. 
\item[c)] \emph{LRBap} if $f$ is LBap and RBap.
\item[d)]  \emph{Bap} if $f\in C(G,X)$ and $P(f,U)$ is weakly rd.
\end{enumerate}
\end{definition}
Here the Bap must not be confused with the more general Besicovitch ap (see e.g. \cite[p. 228/229]{lit38}).

\setcounter{definition}{11}
\begin{lemma}\label{lemma3.12}
For $G,X, \mathscr{U}, f$ as in Lemma 3.10, (3.3) and (3.4) are equivalent, with
\begin{enumerate}
\item[(3.3)] For any $U\in\mathscr{U}$, the $P_L(f,U)$ is left-relatively dense.
\item[(3.4)] For any $U\in\mathscr{U}$, the $P_L(f,U)$ is right-relatively dense.
\end{enumerate}
The same holds for the $P_R$ and the $P$.
\begin{proof}
As for Lemma 3.10, with $K $ also $K^{-1}$ is compact.
\end{proof}
\end{lemma}
Replacing "left rd" in Definition 3.11a by "right rd" gives the same function class, similarly for Definition 3.11b. 

For the $P(f,U)$ see also Corollary 4.3.
\begin{definition}\label{def3.13a}
$f: G\rightarrow X$  is called
\begin{enumerate}
\item[a)]  \emph{left-uniformly  continuous (Luc)} means  to each $U\in \mathscr{U}$ there exists a neighborhood $V$ of the unit $e$ in $G$ such that $$(f(x),f(y)) \in U\text{ for all } x,y\in G \text{ with } xy^{-1}\in V.$$
\item[b)]  \emph{Ruc} if  to each $U\in \mathscr{U}$ there exists a neighborhood $V \in \mathscr{N}(e)$ with $$(f(x),f(y)) \in U\text{ if } x,y\in G \text{ with } x^{-1}y\in V.$$
\item[c)] \emph{LRuc} if  $f$ is Luc and Ruc.
\item[d)]  \emph{uc} if to each $U\in\mathscr{U}$ there is a neighborhood $V$ of $e$ in $G$ such that $$(f(xvy),f(xy))\in U \text{ for all } x,y \in G \text{ and } v\in V.$$
\end{enumerate}
\end{definition}
Such $f$ are $\in C(G,X)$, uniformly continuous $f$ are Luc and Ruc; see also Corollary 4.3 and Question 3, $a\Leftrightarrow b\Leftrightarrow c\Leftrightarrow d$.

However, Luc does in general not imply Ruc: See Example 6.15. 

For abelian $G$, the 4 definitions L,R,LR, - coincide in all the Definitions above.
\section{Classical Characterizations of Almost Periodicity}
\begin{theorem}\label{theorem4.1}
For $G$ topological group, $X$ uniform space and arbitrary $f: G\rightarrow X$, for $f$ the following 25 properties of section 3 are equivalent:

LMap (left Maak almost periodic), RMap, LRMap, Map\\
LNap (left Neumann almost periodic), RNap, LRNap, Nap\\
Lnnap (left net-subnet ap), Rnnap, LRnnap, nnap\\
Lsnap (left sequence-subnet ap), Rsnap, LRsnap, snap\\
LFBap (left F-Bohr ap), RFBap, LRFBap, FBap\\
LBap and Luc (left Bohr ap and left uniformly continuous), RBap and Ruc, LRBap, Bap and uc, Bap.

\end{theorem}
With Lemmas 3.10 and 3.12 one gets further equivalent properties.

For abelian $G$, with Lemma 5.2 the "and Luc" resp. "and Ruc" can be ommitted, else not by Example 6.15.

\begin{definition}\label{def4.2}
An $f:G\rightarrow X$ with (2.1) which has one and therefore all of the 25 properties of Theorem 4.1 will be called \emph{almost periodic (ap)}, the set of all such $f$ will be denoted by \emph{$AP(G,X)$}, with the uniformity induced by that of $C(G,X)$, (2.4), (2.5), and sup-norm if $X$ is normed.
\end{definition}
$AP(\mathbb{R}, \mathbb{C})=\{(L)Bap\}$ gives then Bohr's ap functions \cite{lit14}.

\begin{corollary}\label{cor4.3}
With (2.1) and fixed $f\in C(G,X)$ the following 6 statements are equivalent ("All P are" means "For each $U\in\mathscr{U}$ the $P(f,U)$ is ": see Definitions 3.6, 3.7, 3.8):
\begin{quote}
(1) All P are LFrd (2)  All P are RFrd\\ 
(3) All P are Lrd   (4) All P are Rrd \\
(5) All P are weakly Frd (6) All P are weakly rd.\end{quote}
\end{corollary}

\begin{corollary}\label{cor4.4}
With (2.1), (2.3), (2.4) the $AP(G,X)$ of Definition 4.2 is closed in $C(G,X); \ AP(G,X)$ is complete if and only if $X$ is complete.
\end{corollary}

\begin{corollary}\label{cor4.5}
With (2.1) the $AP(G,X)$ is left and right invariant, i.e. with $ f\in AP(G,X)$ and $a\in G$ also the translates ${}_a f$ and $f_a \in AP (G,X)$ and $\tilde f\in C(G,X), \tilde f(x):=f(x^{-1})$.
\end{corollary}
\begin{corollary}\label{cor4.6}
If $f\in AP(G,X), f(G)$ is totally bounded and $f$ is strongly uniformly continuous, i.e. to each $U\in \mathscr{U}$ exists a neighborhood $W$ of unit $e$ in $G$ with 
\begin{equation}(f(usvtw), f(st))\in U \text{ if } s,t\in G, \  u,v,w \in W.\end{equation}
\end{corollary}
\begin{corollary}\label{cor4.7}
With (2.1) and $Z$ another uniform space, $f\in AP(G,X)$ and $g\in AP(G,Z)$ imply
\begin{equation}(f,g) \in AP(G,X\times Z).\end{equation}
\end{corollary}
\begin{corollary}\label{cor4.8}
If $G,X,Z$ are as in Corollary 4.7 and $T: X\rightarrow Z$ is uniformly continuous, then
\begin{equation}T(AP(G,X)) \subset AP(G,Z). \end{equation}
T uc: To $U\in \mathscr{U}_Z$ exists $V\in \mathscr{U}_X$ with $(T(u), T(v)) \in U$ if $(u,v)\in V$.
\end{corollary}
\begin{corollary}\label{cor4.9}
With (2.1), if $X$ is complete resp. topologically complete resp. a topological group for which the left and right uniformities (Example 6.2) coincide resp. a topological linear space over a topological field resp. a Banachspace, then with pointwise defined operations so is $AP(G,X)$.
\end{corollary}
\begin{corollary}\label{cor4.10}
If $G$ is a topological group and $X$ a topological linear space over a topological field $F$, $\varphi \in AP(G,F)$ and $f\in AP(G,X)$, then $\varphi f \in AP(G,X)$.
\end{corollary}
\begin{corollary}\label{cor4.11}
With (2.1), if to $f\in C(G,X)$ exists $M\subset X$ with $f(G)\subset M$ and $M$ complete and satisfying ($A_0$) of section 2, then for $f$ the following 4 properties (Definitions 3.5)
\begin{equation}\text{ ssap, LRssap, Lssap, Rssap }\end{equation}
are equivalent and equivalent with $f\in AP(G,X)$ (Definition 4.2).
\end{corollary}

So for example, if $X$ is a topologically complete group satisfying ($A_0$) for $\mathscr{U}_L$ [or $\mathscr{U}_R$] (see (6.2)), all 29 properties of Theorem 4.1 and (4.4) are equivalent; this holds especially for $X$ a complete metrizable group, so for Banach spaces.

Sometimes however "complete" is too strong an assumption: A Banach space $B$ in its weak topologiy is topologically complete if (and only if) $B$ is reflexive (Example 6.7); but $(B,w)$ is complete if and only if $B$ is finite-dimensional \cite[p. 55, III.1.(5f)]{lit21} (the same holds for dual $B'$ in its $w'$-topology).

Also the more familiar ($A_1$) axiom is in many cases too strong:
One can show that if the weak topology of a normed space $X$ satisfies ($A_1$), then $X$ is finite-dimensional; the same holds for dual $B'$ in the $w'$-topology, $B$ a Banach space. See however Example 6.7.
\section{The Proofs}
\textbf{Proof of Theorem 4.1}, for fixed $G,X$ and $f\in C(G,X)$:\\

\noindent \textbf{I. ZMap $\Leftrightarrow$ ZNap for $Z\in \{L,R,LR,-\}$ where $Z= - $ means Map $\Leftrightarrow$ Nap:}

This follows from the Definition 3.1 of ZMap, Definition 3.2 of ZNap and of totally bounded (2.6);
for example $f$ Nap $\Rightarrow$ $f$ Map: By (2.2) to $U\in \mathscr{U}$ there is $V\in \mathscr{U}$ with $V=V^{-1}$ and $V\circ V\subset U$, $f$ Nap means there are $a_1, \dots, a_n \in G$ so that to each $u\in G$ there is $j_u \in \{1,\dots,n\}$ with $(f(sa_{j_u}t), f(sut))\in V$ for all $t\in G$. Define $A_j:=\{u\in G \mid (f(sa_jt),f(sut))\in V$ for all $s,t \in G$\}, then $\bigcup_{j=1}^n A_j = G$ and if $u,v \in A_j$ for some $j$, $(f_{[a_j]}, f_{[u]}) \in V, (f_{[a_j]},f_{[v]}) \in V$, on $G\times G$; this implies $(f_{[u]},f_{[a_j]}\in V^{-1}=V, (f_{[a_j]},f_{[v]})\in V$, so $f_{[u]},f_{[v]})\in V \circ V\subset U$ on $G\times G$, ie.e $f$ is Map.\\

\noindent \textbf{II. ZNap $\Leftrightarrow$ Znnap $\Leftrightarrow$ Zsnap for $Z \in \{L,R,LR,-\}$:}

This follows from Theorem 2.1, since the $C(G,X)$ resp. $C(G\times G, X)$ used in Definitions 3.2 are uniform spaces by section 2.\\

For the following we need
\begin{lemma}\label{lemma5.1}
(Maak \cite[p. 143, Satz 4]{lit38}). For $G,X$ as in (2.1), $f$ left Maak ap implies $f$ Maak ap.
\begin{proof}
By (2.2) to $U\in\mathscr{U}$ there exist $V, W\in\mathscr{U}$ with $V\circ V\subset U, W\circ W \subset V$ and $W=W^{-1}$. To $W$ exist $n\in \mathbb{N}$, non-empty $A_1, \dots, A_n\in G$ with $ \bigcup_{j=1}^n A_j=G$ so that $(f(ut),f(vt))\in W$ for all $t\in G$ whenever $u,v \in \text{ same } A_j$. Choose $a_j\in A_j$, then $(f(utd),f(a_jtd))\in W $ if $u\in A_j, \ t,d \in G, 1 \leq j\leq n$. Now the $g_k:= {}_{a_k} f , g_k(t):= f(a_kt), $ are again LMap (the $a_k^{-1}A_j$ satisfy Definition 3.1a for $g_k, W$), with $B_{kj}:=a_k^{-1} A_j$ one has $G=\bigcup_{j=1}^n B_{kj}$ and $(g_k(ut), g_k (vt))\in W$ if $u,v \in B_{kj}$ for the same $j, t\in G, j\in I:=\{1,\dots, n\}, k\in I$.

For $\beta=(\beta_1, \dots, \beta_n)\in I^n$ define $C_{\beta}:= \bigcap_{k=1}^n B_{k{\beta_k}}$.

Then $\bigcup_{\beta\in I^n} C_{\beta}= G:$ To $x\in G$ and $k\in I$ exists $\beta_{k,x}\in I$ with $x\in B_{k\beta_{kx}}$, so $x\in C_{\beta(x)}, \beta(x):=(\beta_{1x},\dots, \beta_{nx})$.
The $C_{\beta}$ form a partition of $G$ in the sense of Definition 3.1a for the $g_1, \dots,g_n$ simultaneously:

If $u,v\in \text{ same } C_{\beta}, \beta=(\beta_1, \dots,\beta_n)$, then $u,v \in B_{k\beta_{k}}$ for all $k\in I$ imply $(g_k(ut),  g_k(vt)) \in W$ for all $t\in G, k\in I$.
These $C_{\beta}, \beta \in I^n$  now also form a partition for $f$ and $U$ in the sense of Definition 3.1d:
Fix $\beta \in I^n, \ s,t \in G, u,v \in \text{ this } C_{\beta}$; to $s$ there is $k_s \in I$ with $s \in A_{k_s}$, so 
\begin{equation}(f(sut),f(a_{k_s} ut))\in W, \ (f(a_{k_s}vt), f(svt))\in W. \end{equation}
Furthermore, since $u,v \in \text{ same } C_{\beta}, (C_{\beta})$ a partition for all $g_k$ by the above, one has
\begin{equation}(f(a_{k_s}ut),f(a_{k_s} vt))\in W.\end{equation}

(5.1), (5.2) and $W\subset W\circ W \subset V$ imply $(f(sut),f(a_{k_s}vt))\in W\circ W \subset V$ and $(f(a_{k_s}vt),f(svt))\in W\circ W \subset V$ and so 
$$(f(sut),f(svt)) \in V\circ V \subset U \text{ for any } s, t \in G.$$
\end{proof}
\end{lemma}

\noindent \textbf{III. LMap $\Leftrightarrow$  RMap $\Leftrightarrow$ LRMap  $\Leftrightarrow$  Map:}

If $\tilde{f} (x):= f(x^{-1}), f $ Rmap implies $\tilde{f}$ LMap (with $A_j^{-1}$), so with Lemma 5.1 one gets:
$f$ LMap  $\Rightarrow$ $f$ Map  $\Rightarrow$ $f$ LRMap  $\Rightarrow$ $f$  RMap $\Rightarrow$ $\tilde f$ LMap  $\Rightarrow$ $\tilde f$ Map  $\Rightarrow$ $\tilde f$ RMap  $\Rightarrow$ $\tilde{\tilde f} $ LMap   $\Rightarrow$ $f=\tilde{\tilde{f}}$ LMap.\\

\noindent \textbf{IV. The $16  \ ZS$, $Z\in \{L,R,LR,-\}, S\in \{Map,Nap, nnap, snap\}$ are all equivalent:}

This follows from I, II, III.\\

\noindent \textbf{V. Map $\Rightarrow$ (For all $U\in\mathscr{U}$ the $P(f,U)$ is left F-relatively dense) $\Rightarrow$ FBap  $\Rightarrow$ Bap:}

By Definition 3.1d to $U\in \mathscr{U}$ there exist $A_1, \dots, A_n \subset G$ and $a_j \in A_j, 1\leq j\leq n$ with $\bigcup_{j=1}^n A_j=G$ and $(f(sut),f(sa_jt))\in U$ for $s,t \in G$ and $u\in A_j$; replacing $s$ by $sa_j^{-1}$ gives $(f(sa_j^{-1}ut), f(st))\in U$ for $s,t \in G$, i.e. $a_j^{-1}u\in P(f,U)$ if $u\in A_j, 1\leq j \leq n$. Now if $x\in G$ there is $j_x$ with $x\in A_{j_x}$, so $a_{j_x}^{-1}x \in P(f,U)$ or $ x\in FP(f,U)$ with $F:=\{a_1,\dots,a_n\}: P(f,U)$ is left F-relatively dense and thus weakly Frd and weakly rd as needed in Definitions 3.9d, 3.11d.\\

\noindent \textbf{VI. All $P(f,U)$ LFrd $\Rightarrow$ LFBap $\Rightarrow$ LMap:}

2nd $\Rightarrow$: With $U,V$ of (2.2) to  $P:=P_L(f,V)$ exists $F=\{a_1, \dots, a_n\}$ with $PF=G$  by Lemma 3.10. Define $A_j=Pa_j$, then $\bigcup_{j=1}^n A_j =G$. If $u\in A_j, t\in G, $ then $u=\tau a_j$ with $\tau \in P$, so $(f(ut), f(a_jt))=(f(\tau a_jt), f(a_jt))\in V$. So if $u,v\in \text{ same } A_j$ , $(f(a_jt), f(vt))\in V^{-1}=V$, or $(f(ut), f(vt))\in V\circ V\subset U$, any $t\in G$, i.e. $f$ is LMap (Definition 3.1a).\\

\noindent \textbf{VII. All $P(f,U)$ LFrd $\Rightarrow$ RFBap $\Rightarrow$ RMap:}

This follows from Lemma 3.10 as in VI, or with VI for $f(x^{-1})$.\\

\noindent \textbf{VIII. All $P(f,U)$ LFRrd $\Leftrightarrow$ LRFBap $\Leftrightarrow$ LFBap $\Leftrightarrow$ RFBap $\Leftrightarrow$ Map:}

This follows from V, VI, VII and III.\\

\noindent \textbf{IX. All $P(f,U)$ LFrd $\Rightarrow$ all $P(f,U)$ Lrd $\Rightarrow$ LRBap $\Rightarrow$ LBap and $f$ Luc:}

The first two "$\Rightarrow$" follow from the definitions and Lemma 3.12; Lemma 5.2 below gives the last "$\Rightarrow$".
\begin{lemma}\label{lemma5.2}
With (2.1), $f$ LBap implies $f$ Ruc, $f$ RBap implies $f$ Luc.
\begin{proof}
For given $ U\in\mathscr{U}$ one has to show $(f(x), f(y))\in U$ if $x^{-1}y\in $ suitable neighborhood $W$ of $e$ in $G$ (Definition 3.13b). By (2.2) to $U$ exist $V, V_1\in \mathscr{U}$ with $V_1=V_1^{-1} , V\circ V\subset V_1, V_1\circ V_1\subset U$. With $P:=P_L(f,V)$ there exists a compact $K\subset G$ with $PK=G$, using Lemma 3.12. Since $f$ is continuous on $G$, to each $a\in K$ exist a neighborhood $N_a$ of $e$ in $G$ with $(f(aw),f(a))\in V$ if $w\in N_a$. To $N_a$ exists an open neighborhood $W_a$ of $e$ in $G$ with $W_aW_a \subset N_a, a \in aW_a, aW_a$ is open. Since $K\subset \bigcup_{a\in K} a W_a$, there exist $a_1, \dots, a_n\in K$ with $K\subset \bigcup_{j=1}^na_j W_{a_j}$. Define $W=\bigcap_{j=1}^n W_{a_j}$, a neighborhood of $e$ in $G$. If $x,y \in G$ with $w:= x^{-1}y \in W $, one has $y=xw$. $x\in G=PK$ gives $x=\tau_x a_x$ with $\tau_x \in P, a_x \in K$, so $a_x \in a_{j_x}W_{a_{j_x}}, x=\tau_x a_{j_x} w_x, w_x \in W_{a_{j_x}}$. Then $(f(x), f(a_{j_x} w_x))=(f(\tau_x a_{j_x} w_x), f(a_{j_x}wx))\in V, \ (f(a_{j_x}w_x), f(a_{j_x}))\in V$ since $w_x\in W_{a_{j_x}} \subset N_{a_{j_x}}, $ so 
$$(f(x),f(a_{j_{x}}))\in V\circ V \subset V_1.$$
$(f(y), f(a_{j_x}w_xw))=(f(\tau_x a_{j_x} w_x w), f(a_{j_x}w_xw))\in V$,
$(f(a_{j_x}w_xw),f(a_{j_x}))\in V$ since $w_x\in W_{a_{j_x}}, w\in W \subset W_{a_{j_x}}, W_{a_{j_x}}W_{a_{j_x}}\subset N_{a_{j_x}}.$ This gives $(f(y),f(a_{j_x}))\in$ $V\circ V \subset V_1$ or $ (f(a_{j_x}), f(y))\in V_1^{-1}=V_1$. Together one gets $(f(x),f(y))\in V_1\circ V_1\subset U$, $f$ is Ruc.\\
The second implication follows by applying the first to $f(x^{-1})$.
\end{proof}
\end{lemma}

\noindent \textbf{X. LBap and $f$ LUC $\Rightarrow$ LMap:}

With (2.2) to $U\in \mathscr{U}$ exist $V\in \mathscr{U}$ with $V=V^{-1}$ and $V\circ V\subset U$, then $V_1\in\mathscr{U}$ with $V_1\circ V_1\subset V$. By Lemma 3.12, to $V_1$ exists a compact $K\subset G$ with $ PK=G, P:=P_L(f,V_1)$. $f$ Luc gives an open neighborhood $W$ of $e$ in $G$ with $(f(x),f(y))\in V_1$ for all $x,y \in G$ with $xy^{-1}\in W$. $K$ compact gives $a_1, \dots, a_n\in K$ with $K\subset \bigcup_{j=1}^n Wa_j$. Define $A_j:=PWa_j, 1\leq k \leq n$.\\ Then $\bigcup_{j=1}^n A_j = G: $ $ x\in G=PK$ gives $x=\tau a$ with $\tau \in P, a\in K$, so $x=\tau w a_j$ for some $j$ and $w\in W$, i.e. $x \in $ this $A_j$.

\begin{sloppypar}
	The $A_1, \dots, A_n $ satisfy Definition 3.1a for $f $ and $U$: Fix $j\in \{1,\dots,n\}$, assume $ u,v\in A_j, t\in G: $ then $u=\tau w a_j$ with $\tau \in P, w\in W$, so $A:=(f(ut), f(wa_jt))=(f(\tau w a_j t), f(w a_j t))\in V_1$, similarly $B:=(f(\tau a_j t), f(a_jt))\in V_1$, $C:= (f(wa_jt), f(a_j t))\in V_1$ by the definition of $W$, $wa_jt(a_jt)^{-1}=w \in W$. $A\in V_1, C\in V_1$ imply $(f(ut),f(a_jt))\in V_1\circ V_1 \subset V$. Since $v$ also $\in A_j$, one gets $(f(vt), f(a_jt))\in V$ therefore $(f(ut),f(vt))\in V\circ V\subset U$.\\
\end{sloppypar}

\noindent \textbf{XI. LRBap $\Rightarrow$ RBap and $f$ Ruc  $\Rightarrow$  RMap:}

Apply IX and X to $f(x^{-1})$ (see III).\\

\noindent \textbf{XII.  All $P(f,U)$ Lrd $\Leftrightarrow$ LBap and $f$ Luc  $\Leftrightarrow$ RBap and $f$ Ruc  $\Leftrightarrow$ Map:}

This follows from IX, X, XI, III, V.\\

\noindent \textbf{XIII.  Bap $\Rightarrow$  FBap and $F$ uc:}

We first need 
\begin{lemma}\label{lemma5.3}
To $f\in C(G,X), \ K \text{ compact }\subset G$, $U\in \mathscr{U}$ there exists a neighborhood $W$ of $e$ in $G$ such that
$$(f(awb),f(ab))\in U \text{ for all } w\in W, a,b \in K.$$
\begin{proof}[Proof by contradiction:] 
Else there exists $U_0\in \mathscr{U}$ such that to each neighborhood $N$ of $e$ there are $a_N, b_N\in K, w_N\in N$  with $(f(a_Nw_Nb_N)f(a_N,b_N))$ not in $U_0$. Here $(a_N)_{N\in \mathscr{N}(e)}, (b_N)$ are nets from $K$,$(w_N)$ is a net from $G$ with $w_N \rightarrow e$; hence there exists subnets $(a_M), (b_M), (w_M)$ (same $M$ for all three) with $a_M \rightarrow a\in K, b_M \rightarrow b\in K, w_M\rightarrow e$. This implies $a_Mw_Mb_M \rightarrow aeb=ab$, so $f(a_Mw_Mb_M) \rightarrow f(ab), a_Mb_M \rightarrow ab, f(a_Mb_M)\rightarrow f(ab)$. Now by (2.2) there exists $U\in \mathscr{U}$ with $U=U^{-1}$ and $U\circ U\subset U_0$. To $U$ exists $M_0$ with $ (f(a_{M_0}w_{M_0}b_{M_0}), f(ab))\in U, (f(a_{M_0}, b_{M_0}), f(ab))\in U$, implying $(f(a_{M_0}w_{M_0}b_{M_0}, f(ab)))\in U\circ U\subset U_0$. Since $M_0$ is a $N\in\mathscr{N}(e)$, this is a contradiction with our choice of the $a_N, w_N, b_N$. 
\end{proof}
\end{lemma}
\begin{lemma}\label{lemma5.4}
$f$ Bap implies $f$ uc (Definition 3.13d).
\begin{proof}
By Definitions 3.11d and 3.7d, to $V\in \mathscr{U}$ there exists a compact $K\subset G$ with $KPK=G, P:=P(f,V), K=K_1\cup K_2$. With (2.2), to $U$ exist $V_j\in \mathscr{U}$ with $V_j=V_J^{-1}, j=1,2,3 $ and $V_1\circ V_1\subset U, V_2\circ V_2 \subset V_1, V_3\circ V_3 \subset V_2$. By Lemma 5.3 to $V_2$ there exists $W\in \mathscr{N}(e)$ sucht that $(f(xwy), f(xy))\in V_2$ for all $w\in W$, all $ x,y \in L:= KK$, again compact (all nets from $L$ have in $L$ convergent subnets).

To $s,t \in G$ exist $a,b,c,d \in K$ and $\tau, \tau ' \in P$ with $s=a\tau b, t=c\tau 'd, $ now $V=V_3, \ P=P(f,V_3)$; then $(f(swt), f(abwt)) \in V_3, (f(abwt),f(abwcd))\in V_3$. 
This gives $(f(swt), f(abwcd))\in V_3\circ V_3 \subset V_2$. If now $w\in $ the above $W$, one gets $(f(abwcd),f(abcd))\in V_2$, so $(f(swt),f(abcd))\in V_2\circ V_2 \subset V_1$; from this $(f(st), f(abcd))\in V_1=V_1^{-1}$ if $w=e$.
So finally $(f(swt), f(st))\in V_1\circ V_1 \subset U$ if $w\in W, s,t \in G$, as desired.
\end{proof}
\end{lemma}

\begin{lemma}\label{lemma5.5}
$f$ Bap implies $f$ FBap and $f$ uc.
\begin{proof}
With Lemma 5.4 we have only to show (Definitions 3.9d and 3.6d) that to each $U\in \mathscr{U}$ there exists a finite $F\subset G$ with $FP(f,U)F=G$. For this, to given $U$ choose $V_1, V_2$ as in the proof of Lemma 5.4. By assumption, to $V_1$ there exists a compact $K\subset G$ with $KP(f,V_1)K=G$. By Lemma 5.4 to $V_2$ there exists $W\in \mathscr{N}(e)$ with $(f(swt), f(st))\in V_2$ if $w\in W, s, t \in G$. Define $P=P(f,V_1), T=WPW$; then if $s,t \in G, \sigma \in T, \sigma=u\tau v$ with $u,v \in W, \tau \in  P$, one has $(f(s\sigma t)f(suvt))=(f(su\tau v t), f(suvt))\in V_1, \ (f(suvt),f(svt))\in V_2, \ (f(svt), f(st))\in V_2$; with $(f(s\sigma t), f(suvt))\in V_1$ this gives $(f(s\sigma t), f(st))\in V_1\circ V_1\subset U$, i.e. $T\subset P(f,U)$.

$T$ is weakly Frd: To the above $K,W$ (the latter can be choosen open) there exist $a_j, b_k\in K$ with $K \subset \bigcup_{j=1}^m a_j W, K\subset \bigcup_{j=1}^n Wb_k$; with these define $F:=\{a_1, \dots, a_m, b_1, \dots, b_n\}$. Then $FTF=G:$ If $x\in G=KPK, x= a_j u \tau v b_k$ with suitable $j,k, \ u,v \in W, \tau \in P, $ so $x\in FWPWF=FTF, FTF=G$.
\end{proof}\end{lemma}
\noindent \textbf{XIV. FBap $\Rightarrow$ Map:}

To $U \in \mathscr{U}$ choose $V=V^{-1} \in \mathscr{U}$ with $V\circ V\subset U$ by (2.2); by assumption there exists $F= \{a_1,\dots, a_n\}\subset G$ with $FPF=G, P:=P(f,V)$. 
Define $A_{jk}=a_j P a_k$; then $\bigcup_{j,k} A_{jk}=G$. The $(A_{jk})$ is a U-partition for $f$: 
If $u\in A_{jk}, u= a_j\tau a_k$ with $\tau \in P$, so $(f(sut), f(sa_ja_k t))=(f(sa_j\tau a_k t), f(sa_ja_kt))\in V$ for all $s,t\in G$; if $v\in $ same $ A_{jk}$, also $(f(svt), f( sa_ja_kt))\in V$, 
so together one gets $(f(sut),f(svt))\in V\circ V\subset U$, $f$ is Map.\\

\noindent \textbf{XV. Bap $\Leftrightarrow$ FBap $\Leftrightarrow$ FBap and $f$ uc  $\Leftrightarrow$ Map:}

This follows from XIII, XIV, V.\\

The proof of Theorem 4.1 is now complete.

\begin{proof}[Proof of Corollary 4.3]
With Definitions 3.6, 3.7, 3.8, Lemmas 3.10, 3.12, "$f$ Bap $\Rightarrow$ $f$ Map" of Theorem 4.1 and "$f$ Map $\Rightarrow$ all $P$ are LFrd" of Part V of the proof of Theorem 4.1
\end{proof}

\begin{proof}[Proof of Corollary 4.4]
"AP closed": Definitions 3.1 for $V,U$ with $(V\circ V)\circ V \subset U$ of (2.2). "only if": constants $X\subset AP(G,X)$. \\ "if": If for a Cauchy net $(f_j)$ the $f_j(x) \rightarrow : f(x)$ for all $x\in G$, then $f_j \rightarrow f$ uniformly on $G$, so $f\in C(G,X)$.
\end{proof}

\begin{proof}[Proof of Corollary 4.5]
With $a^{-1}A_j$ resp. $A_j a^{-1}$ resp. $A_j^{-1}$ in Definition 3.1a.
\end{proof}

\begin{proof}[Proof of Corollary 4.6]
$f(G)$ totally bounded follows with Definition 3.2a or (5.4). By Theorem 4.1 the $f$ is uc (Definition 3.13d), repeated use of this and (2.2) give (4.1).
\end{proof}

\begin{proof}[Proof of Corollary 4.7]
With Definition 3.3a, the definition of the uniformity for $X\times X$ \cite[p. 182]{lit31}, and repeated subnets.
\end{proof}

\begin{proof}[Proof of Corollary 4.8]
Obvious.
\end{proof}

\begin{proof}[Proof of Corollary 4.9]
We omit the proof for $X=$ non-abelian group. If $X$ is an abelian topological group, $f,g \in AP(G,X), (s_j)$ a net from $G$, there exists a subnet $(t_{\beta})$ with $(f_{t_{\beta}})$ and $(g_{t_{\beta}})$ Cauchy, commutativity of "$+$" and (2.2) give $(f+g)_{t_{\beta}}$ Cauchy; similar in the other cases. One could also use Corollary 4.7.
\end{proof}

\begin{proof}[Proof of Corollary 4.10]
With Corollary 4.7, the $(\lambda, x)$-commutativity of $\lambda x$ in $(0,0)$ and the total boundedness of $\varphi(G), f(G)$.
\end{proof}
\begin{proof}[Proof of Corollary 4.11]
For $f\in C(G,X)$, with Definition 3.5 and Theorem 4.1 one has
\begin{equation}\text{ssap} \Rightarrow \text{LRssap} \Rightarrow \genfrac{\{}{\}}{0pt}{1}{\text{Lssap} \Rightarrow \text{Lsnap}}{\text{Rssap} \Rightarrow \text{Rsnap}} \Rightarrow  \text{snap}.\end{equation}
It remains to show  snap $\Rightarrow$ ssap.
By Theorem 4.1, $f$ is RNap, i.e $\{f_t \mid t\in G\}$ is totally bounded in $C(G,X)$; this and (2.4), (2.5) imply
\begin{equation} f\in AP(G,X) \Rightarrow f(G) \text{ totally bounded in }X.\end{equation}
Then the closure $\overline{f(G)}$ is also totally bounded with (2.2), $K:=\overline{f(G)} \subset \overline M =M $ of Corollary 4.11, $K$ is also complete. By 32 Theorem of \cite[p. 198]{lit31} $K$ is compact.
Since $M$ satisfies $(A_0)$, $K$ also satisfies $(A_0)$. By Lemma 5.6 below, $K$ satisfies $(A_1)$, there exists a sequence $(U_n)_{n\in\mathbb{N}}$ of $U_n \in\mathscr{U}$ which form a basis for $\mathscr{U}\big|_{K}$, one can assume $U_{n+1} \subset U_n$ , $n\in \mathbb{N}$.\\
Let there now be given arbitrary sequence $(s_n)_{n\in\mathbb{N}}$ from $G$. $f$ snap means there exists a subnet $(s_{n(j)})_{j\in J}$ so that $(f_{[n(j)]})$ is Cauchy in $C(G\times G, X)$. $(s_{n(j)})$ subnet of $(s_n)$ means
\begin{equation}\text{to any } m\in\mathbb{N} \text{ exists } j_m\in J \text{ with } m\leq n(j) \text{ if }j_m\leq j.\end{equation}
$(f_{[n_j]})$ Cauchy implies: If $m\in\mathbb{N}$ is given, and, by (2.2), to $U_m$ one chooses $V_m \in\mathscr{U}$ with $V_m\circ V_m \subset U_m$, then 
\begin{multline}\text{to } V_m  \text{ exists } k_m\in J, j_m\leq k_m\text{ , with }(f_{[s_{n(u)}]}(s,t), f_{[s_{n(v)}]} (s,t) \in V_m\\ \text{ for all }s,t\in G\text{  if }u,v \in G\text{  with }k_m \leq u\text{  and }k_m\leq v.\end{multline}
Define (axiom of choice) $q(m):= n(k_m), m\in\mathbb{N}$. Then $(s_{q_m})_{m\in\mathbb{N}}$ is a subsequence of $(s_n)$: If, for given $n$, $n\leq m $, then $j_m\leq k_m$ and (5.5) imply $n\leq m \leq n(k_m)=q(m)$.
For fixed $s,t\in G, \ (f(ss_{n(j)} t))_{j\in J}$ is Cauchy in $K=\overline{f(G)}; K$ being complete, there exists $F(s,t)\in K$ with $f(ss_{n(j)}t) \rightarrow F(s,t)$. So to $m$ there is $j=j(s,t,m) \in J$ with $(f(ss_{n(j)}t), F(s,t))\in V_m$ and $k_m\leq j$. With (5.6) one gets $(f(ss_{q(m)}t), F(s,t))\in V_m\circ V_m \subset U_m$, for all $s,t \in G, m\in \mathbb{N}$.

So if $p\in \mathbb{N}$ is given, $p\leq m \in \mathbb{N}, (f(ss_{q(m)}t), F(s,t))\in U_m \subset U_p$, for all $s,t\in G$. This means, since the $(U_p)_{p\in\mathbb{N}}$ form a basis for $U\big|_{K}$, that $(f_{[s_{q(m)}]})_{m\in\mathbb{N}}$ is Cauchy in $C(G\times G, X)$.
So $f$ snap does indeed imply $f$ ssap.
\end{proof}

\begin{lemma}\label{lemma5.6}
If $(X, \mathscr{U})$ is a compact uniform space which satisfies $(A_0)$, then $(X, \mathscr{U})$ satisfies $(A_1)$\end{lemma}
See also \cite[Lemma 11, and the remarks after it]{lit26}.
\begin{proof}
By $(A_0)$ there exists a sequence $(U_n)_{n\in\mathbb{N}}$ from $\mathscr{U}$ with $\bigcap_{n=1}^{\infty} U_n = \bigcap \mathscr{U}$. The $(U_n)$ form a basis for $\mathscr{U}$: Else there exist $U_0\in\mathscr{U}$ and $(u_n, v_n)\in U_n$ with $(u_n, v_n)$ not in $U_0, n\in\mathbb{N}$. 
Since $X$ is compact, there exists a subnet $(u_{n(j)})_{j\in J})$ of $(u_n)$ and $u\in X$ with $u_{n(j)} \rightarrow u$. Similarly there are a subnet $(v_{m(k)}):=(v_{n(j(k))})$ of $(v_{n(j)})$ and $v\in X$ with $v_{m(k)}\rightarrow v$.
Since $(m(k))$ is a subnet of $(n(j))$, $ u_{m(k)} \rightarrow u$. 
Now given $U_n$, choose $V_n$ with (2.2); by the above there is $k$ with $(u,u_{m(k)})\in V_n$ and $(v,v_{m(k)})\in V_n=V_n^{-1}$, so $(u,v) \in V_n \circ V_n \subset U_n$. Since this holds for any $n$, $(u,v) \in \bigcap_{n=1}^{\infty}=\bigcap \mathscr{U}$. This implies, with (2.2), $v_{m(k)} \rightarrow u$, and then $(u_{m(k)}, v_{m(k)})\in U_0$ for some $k$, a contradiction.
\end{proof}
\section{Examples}
\noindent 1. Any semimetric space $(X,d), 0=d(x,x)\leq d(x,y)=d(y,x)\leq d(x,z)+d(z,x) \in \mathbb{R},\ x,y,z \in X$, becomes a uniform space with $$\mathscr{U}_d = \{U\subset X\times X \mid U \supset \{(x,y) \mid x,y \in X, d(x,y)<\epsilon\} \text{ for some } \epsilon \in (0, \infty)\},$$ it satisfies $(A_1)$ and so $(A_0)$.

 $C(G,X)$ of (2.3) is then again semimetric with $d(f,g) : = \sup\{d(f(s),g(s)) \mid s \in G \}$ ; the same holds for $X$ seminormed linear with  $||f|| : =
        \sup \{ ||f(s)||_X :  s \in G \} $
\\

\noindent 2. Any topological group $G$ becomes a uniform space with $\mathscr{U}_L$, and also with $\mathscr{U}_R$, where $\mathscr{U}_L:=\{U\subset G\times G \mid U \supset {}_V U \text{ for some } V\in \mathscr{N}(e)\}$ with $ \mathscr{N}(e)=$ set of neighborhoods of unit $e$ of $G$, and 
$$_V U := \{(x,y)\in G\times G \mid x^{-1} y \in V\}, \  V\in \mathscr{N}(e),$$
similarly $\mathscr{U}_R$ with $U_V:=\{(x,y)\in G\times G \mid xy^{-1} \in V\}$.
Both uniformities induce the same topology in G, the original one. They coincide if and only if \cite[p. 21-23]{lit30}
\begin{equation}\text{to any }V\in \mathscr{N}(e) \text{ there is }W\in \mathscr{N}(e)\text{  with }xWx^{-1}\subset V \text{  for all }x\in G.\end{equation}
For non-abelian non-compact $G$ however, $\mathscr{U}_L $ and $\mathscr{U}_R$ are in general not equivalent: See example 16.\\
Special cases: $X$ a linear topological space over a topological field, for example $X$ seminormed linear, a Banach space $B$, or distribution spaces as $\mathscr{D}(\Omega, \mathbb{C}), \mathscr{D}'(\Omega, B)$, $\mathscr{S}(\mathbb{R}^n, \mathbb{C}), \mathscr{S}'(\mathbb{R}^n, B)$, $\Omega \text{ open } \subset \mathbb{R}^n$.\\

\noindent 3. If $X$ is a topologically complete topological group satisfying $(A_0)$ for $\mathscr{U}_L$ (or $\mathscr{U}_R$), then for this $X$ and any $G$ the 29 conditions of Theorem 4.1 and Corollary 4.11 are equivalent; this holds especially for $X$ complete metric group, so any Banach space.\\
If $X$ is a topological abelian group, the $X$ with the induced uniformity (see example 2) satisfies $(A_0)$ if and only if 
\begin{equation}\text{  there exists a sequence of }V_n\in \mathscr{N}(e)\text{  with }\bigcap_{n=1}^{\infty} V_n =\bigcap \mathscr{N}(e). \end{equation}

\noindent 4. If $G$ is a topological group, $V\in \{\mathbb{R},\mathbb{C}\}, f:G\rightarrow V$, define $F:G\rightarrow V^G$ by $(F(s))(t)=f_s(t)=f(ts), \ s,t \in G$. Then the following 3 statements are equivalent, with $B=AP(G,V), BUC(G,V) $ (bounded uniformly continuous functions), $BC(G,V)$ (bounded continuous functions), or $l^{\infty}(G,V)$, Banach spaces with the sup-norm:

a) $f\in AP(G,V)$, b) $F \in AP(G,B)$, c) $F$ weakly ap, i.e $F\in AP(G,(B,w))$,\\  where $(B,w)$ denotes $B$ in its weak topology \cite[p. 235]{lit32}.\\
( a) $\Rightarrow$ b) $\Rightarrow$ c) follows from the definitions and Corollaries 4.6, 4.8, c) $\Rightarrow$ a) with $\varphi(h):=h(e), h\in B: \varphi \in \text{ dual } B'$.)\\

\noindent 5.  If $X$ is a locally convex real or complex separable linear space which satisfies $(A_0)$ of (6.2), then the weak topology of $X$ satisfies $(A_0)$ also.\\
The proof of Lemma 8 of \cite[p. 203]{lit26} works also for the more general $(A_0)$ here, no separation (Hausdorff) axioms for $X$ are needed.\\
Special case: $X$ seminormed linear separable implies $(X,w)$ has $(A_0)$. \\
So for $f\in C(G,X)$ with $X$ seminormed linear and $f(G)$ separable, e.g. $G$ separable, one has $(A_0)$ for the weak closure of $f(G)$ in $(X, \text{ weak topology})$.\\

\noindent 6. A converse of example 5 is false: $AP(\mathbb{R},\mathbb{C}), BUC(\mathbb{R},\mathbb{C}) , BC(\mathbb{R},\mathbb{C}) $ (see example 4) all satisfy $(A_0)$ in the weak topology, none is separable (=weakly separable).
(\cite[Lemma 1, p.337]{lit28}; AP not separable: if $(f_n)$ dense, choose $\omega$ not in $\bigcup_{n} \text{Bohr-spectrum} (f_n), f_{n_k} \rightarrow e^{i\omega t}$, then $0=\text{mean} (f_{n_k}e^{-i\omega t}) \rightarrow 1.$)\\

\noindent 7. For reflexive Banach spaces $B$, Corollary 4.11 can be applied to any $f\in AP(G,(B,w))$ with  $w$ the weak topology and with $f(G)$ separable.\\
 ($B$ reflexive is equivalent with $B$ weakly topologically complete \cite[p. 38 \mbox{Theorem 1}, p.70 Theorem 2(B)]{lit21}; $(A_0)$ with example 5;\\
here the in $\overline{f(G)}$  induced weak uniformity satisfies even $(A_1)$.)\\

\noindent 8. More information on $AP(G,(B,w))$ with $B$ Banach space can be found in \cite[Theorem 3, corollary I p. 201, Corollary and Lemma 10 p. 204]{lit26}, \cite{lit27}. For example, with all $f\in AP(G,(B,w))$, one has \begin{align*}\{f(G) \text{ separable}\}&=\{ \text{the weak Bohr spectrum of f is countable}\} \\ &\subset \{f \text{ is Lssap}\} \subset \{f \text{ ap}\},\end{align*} with examples showing that the two "$\subset$" are strict (\cite[p. 124, (86), p. 138-141]{lit27}); \\
$AP(G,(B,w))=\{f: G\rightarrow B: y \circ f \in AP(G,F) \text{ for all } y\in B'\},\ F=\mathbb{R} \text{ or } \mathbb{C}$ with Corollary 4.9. ($(B,w)$ complete $\Rightarrow$ dim$B < \infty$ \cite[p.55, (5f)]{lit21}.)
By \cite[(18)]{lit27} for $B$ Banach space the "the $(B,w)$ satisfies $(A_0)$" is equivalent with "dual $B'$ is separable in the w'-topology". See also examples 10 and 14, and \cite[p.200, 203/204]{lit26}. \\

\noindent 9. For $\Omega$ open $\subset \mathbb{R}^n, n\in \mathbb{N}, X$ Banach space, the distribution spaces $\mathscr{D}(\Omega, \mathbb{C})$, $ \mathscr{D}'(\Omega, X)$, $\mathscr{S}(\mathbb{R}^n,\mathbb{C}), \mathscr{S}'(\mathbb{R}^n,X)$ are complete and satisfy $(A_0)$, so Corollary 4.11 is applicable to any $f\in AP(G,T)$, $T$ any of the above distribution spaces, any $G$;
with the exception of $\mathscr{S}(\mathbb{R}^n,X)$ they do not satisfy the first countability axiom $(A_1)$.\\ 
(If $\Omega=\mathbb{R}^n $ and $X=\mathbb{C}$, for "complete" see \cite[p. 66, 71, 235, 238]{lit46}; $(A_0)$ for $\mathscr{D}'(\Omega, X), \mathscr{S}'(\mathbb{R}^n,X)$ follows with $\mathscr{D}(\Omega, \mathbb{C}) $ separable, dense in $\mathscr{S}$ ($\Omega=\mathbb{R}^n)$ and $\mathscr{D}'$ \cite[p. 75]{lit46}; $(A_0)$ for $\mathscr{D}$ follows with $\mathscr{D},\mathscr{D}'$ reflexive and $\mathscr{D}'(\Omega, \mathbb{C})$ separable.)\\
The above holds also for various other Distribution spaces, e.g. $(\mathscr{D}_{L^1})'=\mathscr{D}'_{L^{\infty}}$ is complete and satisfies $(A_0)$ (for definitions see \cite[p. 367 Prop 2.9]{lit7}).\\

\noindent 10. Already for $G=(\mathbb{R} ,+)$ there exist $X$ and $f\in AP(\mathbb{R},X)$ with $f$ not ssap (Definition 3.5, = Bochner's characterisation), $X$ is even a locally convex topological linear space which is topologically complete: $X=l^{\infty}(\mathbb{R},\mathbb{C})$ with the w'-topology induced by $l^1(\mathbb{R},\mathbb{C})$ in $l^{\infty}=\text{ dual } l^1(\mathbb{R},\mathbb{C})', (f(t))(s):=e^{ist}, \ s,t \in\mathbb{R}$.\\
The same holds for this $f$ and $X=\mathbb{C}^\mathbb{R}$ with the topology of pointwise convergence on $\mathbb{R}$, this $X$ is even complete. \\
For such an $f$ with $X=$ Banach space in its weak topology see \cite[p. 128 Example 5]{lit27}.\\
$((l^{\infty}, w')$ top. complete e.g. with the w'-compactness of the unit ball \cite[p. 47 Theorem 3]{lit21} and w'-totally bounded $\Rightarrow$ norm bounded \cite[p. 38 Theorem 1]{lit21}. The $f$ is ap since $f(\cdot)(s) \in AP(\mathbb{R},\mathbb{C}) $ for each fixed $s\in\mathbb{R}$, so $(y, f(\cdot))$ converges uniformly absolutely if $y\in l^1, (y, f(\cdot)) \in AP(\mathbb{R},\mathbb{C})$ by Corollary 4.9 and 4.4.\\
 $f$ not ssap: If $a_n \rightarrow \infty $ and $f(a_n) \rightarrow \varphi$ pointwise on $\mathbb{R}$, then $\varphi$ is a character on $\mathbb{R}$ and Lebesgue-measurable, then continuous, so $\varphi(s)=e^{i\omega s}$ with some $\omega \in \mathbb{R}$ (e.g. \cite[p. 69 Satz 2]{lit38}), with Lemma 9 of \cite[p. 140]{lit27} one gets a contradiction.)\\

\noindent 11. $X=\mathbb{C}^{\mathbb{R}}$ in example 10 shows that Eberlein's theorem (e.g. \cite[p. 58 Theorem 1]{lit21}) becomes false already for complete locally convex topological linear spaces.\\

\noindent 12. Example 10 shows, that Corollary 4.11 becomes false without "$\overline{F(G)}$ has $(A_0)$", for $G, X, f$ of Example 10 even "$F(G)$ has $(A_0)$" holds, one cannot omit the closure of $f(G)$.
$(\varphi, \psi \in f(G), \varphi(r)=\psi(r) $ for rational $r, \varphi,\psi$ continuous $\Rightarrow \varphi=\psi$.)\\

\noindent 13. Also $\overline{f(G)}$ complete is essential for Corollary 4.11:
For $G,f$ of example 10 and $X=BC(G,\mathbb{C})=\{f:\mathbb{R}\rightarrow\mathbb{C} \mid f \text{ continuous bounded}\}$ with the topology of pointwise convergence one has $(A_0)$ even for $X, f\in AP(\mathbb{R},X)$, still $f$ is not ssap by the proof in example 10.\\

\noindent 14. $G,f,X$ of example 10 give simultaneously an explicit example of an almost periodic $f$ with Bohr spectrum $\sigma_B (f)$ uncountable, $=\mathbb{R}$; for other such examples see \cite[p. 140, example 6]{lit27}. (For $\sigma_B$ to be defined a suitable mean $M: AP(G,X)\rightarrow X$ is needed, existing for locally convex topologically complete complex linear spaces \cite{lit11}, \cite{lit37}, \cite{lit39}, \cite{lit23}, \cite[p. 112]{lit27}; with this, $G=\mathbb{R}$, Fourier coefficient $c_{\omega}(f):=M(f\overline{\varphi_{\omega}}), \omega \in \mathbb{R}, \varphi_{\omega}(t)=e^{i\omega t}, \sigma_B(f):=\{\omega \in \mathbb{R} \mid c_{\omega}(f)\ne 0\}$; here $c_{\omega}(f)=\delta_{\omega}, \delta_{\omega}(\omega)=1, \text{ else } 0$.)\\

\noindent 15. In "$f$ L-Bohr ap and L-uniformly continuous" resp. "$f$ RBap and Ruc" in Theorem 4.1 the "L-Uniformly continous" resp. "Ruc" cannot be omitted or replaced by "Ruc" resp. "Luc": 
There exist topological groups $G$ and even $f:G\rightarrow \mathbb{C}$ with $f$ R-Bohr ap, but $f$ not R-uniformly continuous and thus $f$ not L-Bohr ap (Lemma 5.2), $f$ not ap in the sense of Definition 4.2 \cite[Example 1 p. 490, and the references there]{lit41}; a first such example has been given by Wu \cite{lit47}:\\
$G=$ group of Euclidian motions of the plane $\mathbb{R}^2=\mathbb{C}$,\\
$G= \{(z,w) \mid z,w \in \mathbb{C}, |w|=1\}$, with 
\begin{equation}(z',w')(z,w):=(z'+w'z, w'w)\end{equation}
and the usual induced $\mathbb{C}^2$ topology;
\begin{equation}f(z,w):=e^{i\text{Re}(z/w)}, \ (z,w)\in G.\end{equation}
$G$ is even solvable and a Lie group, $f\in C(G,\mathbb{C})$;\\
$$f((z,w)(2\pi k + iy, 1))=f(z+w(2\pi k + iy),w)=f(z,w)$$
 shows $(2\pi \mathbb{Z}+ i\mathbb{R}, 1)\subset P_R(f,\epsilon)$ for any $\epsilon >0$. With $K:= \{(x,w) \mid 0 \leq x \leq 2\pi, w\in \mathbb{C} \text{ with } |w|=1\}$ compact $\subset G$ one has $P_R(f,\epsilon)K=G$ so $f$ is RBap.\\
$f$ is not Ruc: 
\begin{align*}|f((u,v)(z,w))-f(u,v)|&=|(e^{i\text{Re}(u+vz/vw)}-e^{i\text{Re}(u/v)}|\\ &=
|(e^{i\text{Re}(z/w)}-1)e^{i\text{Re}(u/vw)} + (e^{i\text{Re}(u/vw)}-e^{i\text{Re}(u/v)})|\\&=: |A+B|\end{align*}
should tend to $0$ as $(z,w) \rightarrow e=(0,1)$, uniformly in $(u,v)\in G$, especially for $v=1, u\in \mathbb{R}$.\\
Now $A\rightarrow 0$ as $z\rightarrow 0, w\rightarrow 1$, uniformly in $u$ if $v=1$;\\
however, if $u=\pi / (1-w), w\ne 1, |B|=|e^{i\text{Re}((u/w)-u)}-1|=|e^{i\pi / w}-1| \rightarrow 2 \ne 0$: $f$ is indeed not Ruc.\\

\noindent 16. The $G$ of example 15 shows that for non-abelian non-compact groups in general the left uniformity $\mathscr{U}_L$ and the right uniformity $\mathscr{U}_R$ are different, already for metric complete locally compact solvable (="almost abelian") groups.\\

\noindent 17. If $\mathbb{R}_d:=(\mathbb{R}, +)$ in the discrete topology, $AP(\mathbb{R},\mathbb{C})$ strictly $\subset AP(\mathbb{R}_d,\mathbb{C})$, there exist bounded characters $\in AP(\mathbb{R}_d,\mathbb{C})$ which are not even Lebesgue measurable \cite[p. 90/91, Satz 6, Satz 7]{lit38}.


\section{Questions}
\noindent 1. Do there exist $G, X, f: G  \rightarrow X$ so that $P_L(f, U)$ is weakly Frd resp. weakly rd for all $U\in \mathscr{U}$, but $f$ is not almost periodic?\\

\noindent 2. Exist $G$ and $M, K  \subset G$ with $KMK=G$ and $K$ compact (resp. finite), but for no compact (resp. finite) $L\subset G$ one has $LM=G$?\\

\noindent 3.  Do there exist $f:G\rightarrow X$ with $f$ Luc and Ruc, but $f$ not uc (Definitions 3.13)?\\
		(For any fixed topological group $G$, the following 4 properties can be shown to be equivalent:
		\begin{enumerate}
		\item[(a)] For all uniform $X$ and $ f: G\rightarrow X$, $f$ Luc $\Rightarrow f$ Ruc.
		\item[(b)] $G$ satisfies (6.1)
		\item[(c)] Left and right uniformities of $G$ conincide, $\mathscr{U}_L=\mathscr{U}_R$.
		\item[(d)] For all $X$ and $f:G\rightarrow X$, one has: $f$ Luc $\Leftrightarrow$ $f$ Ruc $\Leftrightarrow$ $f$ uc.)
		\end{enumerate}

\noindent 4. Do there exist groups $G, X$ and $f, g \in AP(G,X)$ with $fg$ not ap (see Corollary 4.9)?\\

\noindent 5.  Do there exist $f: G \rightarrow X$ so that to each $U\in \mathscr{U}$ there is a totally bounded $T\subset G$ with $TP_L(f, U)=G$, but $f$ is not ap, though it is (left) uniformly continuous?\\

\noindent 6. Can one generalize Corollary 4.11, e.g. is it true if $\overline{f(G)}$ is only sequentially complete? (See \cite[p. 204]{lit26}, \cite[(86)]{lit27}, Example 6.8)\\

\noindent 7. Under what assumptions are Lssap, Rssap, LRssap, ssap equivalent (not necessarily equivalent with ap)? (Analog to Maak's Lemma 5.1, Lssap $\Rightarrow$ ssap.)\\

\noindent 8. Do there exist $f$ as in examples 6.10, 6.12, 6.14 with $X=$ Hilbert space in its weak topology (and $G$ abelian)?\\

\noindent 9. What of (86) of \cite{lit27} holds if $(B,w)$ is replaced by $(B',w')$? (see example 6.8 and \cite[p. 135, 136, 144, 146]{lit27}.)\\

\noindent 10. Does the weak closure of $f(\mathbb{R})$ in $l^{\infty}(\mathbb{R}, \mathbb{C})$ of example 10 satisfy $(A_0)$?\\

\noindent 11. Does an analog to Lemma 5.1 hold for Maak-ergodic $f:G \rightarrow X$ (See \cite[p. 36, 38]{lit39} for $X=\mathbb{C}$)?\\

\noindent 12. A streamlining of the proof of Theorem 4.1 would be nice.

\noindent 
Mathematisches Seminar\\
Universität Kiel\\
D 24098 Kiel\\
Germany\\
guenzler@math.uni-kiel.de
\end{document}